\documentclass[12pt]{amsart} \textwidth=14.5cm \oddsidemargin=1cm
\usepackage{setspace}
\usepackage{amsthm}
\usepackage[rgb]{xcolor}
\usepackage{times}
\usepackage[linktocpage=true]{hyperref}
\usepackage{mathrsfs}  

\usepackage{xypic}
\usepackage{tikz}
\usepackage{verbatim}
\usepackage{amscd}
\usepackage{amsmath}
\usepackage{amsfonts}
\usepackage{amsxtra}
\usepackage{amssymb}
\usepackage{eucal}
\usepackage{latexsym,todonotes}
\usepackage{stmaryrd}
\usepackage{graphicx}%
\usepackage{dynkin-diagrams}

\usepackage{tikz}
\usepackage[all]{xy}

\theoremstyle{plain}
   \newtheorem{theorem}{Theorem}[section]

 \theoremstyle{definition}
   \newtheorem{example}[theorem]{Example}
\newtheorem{remark}[theorem]{Remark}

\newcommand{\arxiv}[1]{\href{http://arxiv.org/abs/#1}{\tt arXiv:\nolinkurl{#1}}}

\newcommand{\nc}{\newcommand}
 \nc{\A}{\mathcal A}
\nc{\Ainv}{\A^{\rm inv}}
\nc{\aA}{{}_\A}
\nc{\aAp}{{}_\A'}
\nc{\aff}{{}_\A\f}
\nc{\aL}{{}_\A L}
\nc{\aM}{{}_\A M}
\nc{\BH}{{\mathbb H}}
\nc{\Bin}{B_i^{(n)}}
\nc{\bK}{\mathbb K}
\nc{\cp}{\mathcal P}
\nc{\dL}{{}^\omega L}
\nc{\Z}{{\mathbb Z}}
 \nc{\N}{{\mathbb N}}
 \nc{\fZ}{{\mf Z}}
 \nc{\Id}{\text{Id}}
 \nc{\F}{{\mf F}}
 \nc{\Q}{\mathbb{Q}}
 \nc{\la}{\lambda}
 \nc{\ep}{\epsilon}
 \nc{\h}{\mathfrak h}
 \nc{\He}{\bold{H}}
 \nc{\htt}{\text{tr }}
\nc{\iba}{\psi_{\imath}}
 \nc{\n}{\mf n}
 \nc{\g}{{\mathfrak g}}
 \nc{\DG}{\widetilde{\mathfrak g}}
 \nc{\SG}{\breve{\mathfrak g}}
 \nc{\is}{{\mathbf i}}
 \nc{\M}{\mathbb M}
\nc{\V}{\mathbb V}
 \nc{\bi}{\bibitem}
 \nc{\E}{\mc E}
 \nc{\ba}{\tilde{\pa}}
 \nc{\half}{\frac{1}{2}}
 \nc{\hgt}{\text{ht}}
 \nc{\Hom}{\text{Hom}}
 \nc{\res}{\text{res}}
 \nc{\ka}{\kappa}
\nc{\K}{\mathbb{K}}
 \nc{\mc}{\mathcal}
 \nc{\mf}{\mathfrak} 
 \nc{\hf}{\frac{1}{2}}
\nc{\ov}{\overline}
\nc{\ul}{\underline}
\nc{\xx}{{\mf x}}
\nc{\one}{\textbf{1}}
\nc{\Qq}{\Q(v)}

\newcommand{\tk}{\widetilde{k}}

\newcommand{\tMHk}{{\widetilde{\mathcal H}(\mathbb F_q Q,\tau)}}
\newcommand{\va}{\varsigma}
\newcommand{\vs}{\varsigma}
\newcommand{\qbinom}[2]{\begin{bmatrix} #1\\#2 \end{bmatrix} }
\nc{\ua}{\mf{u}}
\nc{\inv}{\theta}
\nc{\mA}{\mathcal{A}}

\newcommand{\I}{\mathbb I}
\newcommand{\II}{\mathcal I}
\newcommand{\Ibw}{\II_{r|m|r}}
\newcommand{\Ianti}{\II_{r|m|r}^{d,-}}
\newcommand{\Iwl}{\II_{\circ}^{-}}
\newcommand{\Iwr}{\II_{\circ}^{+}}
\newcommand{\Ib}{\II_{\bullet}}

\newcommand{\Hy}{\mathscr{H}}

\newcommand{\HB}{\mathscr H_{B_d}}

\nc{\bU}{\textbf{U}}
\nc{\tU}{\widetilde{\textbf{U}}}
\nc{\Udot}{\dot{\bU}}
\nc{\f}{\textbf{f}}
\nc{\fprime}{{}'\textbf{f}}
\nc{\B}{\textbf{B}}
\nc{\Bdot}{\dot{\B}}
\nc{\Dupsilon}{\Upsilon^{\vartriangle}}

\nc{\ipsi}{\psi_{\imath}}
\nc{\Ui}{{\textbf U}^{\imath}}
\nc{\tUi}{\widetilde{\textbf U}^{\imath}}
\nc{\tUiD}{{}^{\mathrm{Dr}}\tUi}
\nc{\Uidot}{\dot{\textbf{U}}^{\imath}}
\nc{\Ktilde}{\tilde{K}}
\nc{\bktilde}{\widetilde{k}}
\nc{\Yi}{Y^{w_0}}
\nc{\bunlambda}{\Lambda^\imath}
\nc{\Iwhite}{\I_{\circ}}
\nc{\ile}{\le_\imath}
\nc{\il}{<_{\imath}}

\nc{\bbW}{{\boldsymbol{W}}}
\nc{\bs}{\mathbf{r}}
\nc{\bwi}{w_{\bullet,i}}
\nc{\tfX}{\widetilde{\Upsilon}}
\nc{\tcT}{\widetilde{\mathcal T}}
\nc{\tT}{\widetilde{T}}
\nc{\tTT}{\widetilde{\mathbf{T}}}

\nc{\qq}{(v_i^{-1}-v_i)}
\nc{\qqq}{(1-v_i^{-2})^{-1}}
\nc{\qqqj}{(1-v_j^{-2})^{-1}}

\nc{\etab}{\eta^{\bullet}}
\newcommand{\Iblack}{\I_{\bullet}}
\newcommand{\wb}{w_\bullet}

\newcommand{\purple}[1]{{\color{purple}#1}}

\numberwithin{equation}{section}

\begin{document}


\title[Quantum symmetric pairs]{Quantum symmetric pairs}

\author[Weiqiang Wang]{Weiqiang Wang}
	\address{Department of Mathematics, University of Virginia, Charlottesville, VA 22904}
	\email{ww9c@virginia.edu}

	


\begin{abstract}
This is a survey of some recent progress on quantum symmetric pairs and applications. The topics include quasi K-matrices, $\imath$Schur duality, canonical bases, super Kazhdan-Lusztig theory, $\imath$Hall algebras, current presentations for affine $\imath$quantum groups, and braid group actions. 
\end{abstract}

\maketitle


\section*{Introduction}

\subsection{Quantum groups}

According to Drinfeld and Jimbo, a quantum group $\bU = \bU_v(\mathfrak g)$ is the quantum deformation (as Hopf algebras) of the universal enveloping algebra $U(\mathfrak g)$, for a semisimple or Kac-Moody Lie algebra $\mathfrak g$. Since their inception in 1985, quantum groups have found numerous applications to diverse areas including mathematical physics, representation theory, algebraic combinatorics, and low dimensional topology. 

We offer a (personal) Top Ten List of highlights in quantum groups, viewed as a part of Lie theory, as follows:

(1) Definition \cite{Dr87, Jim86} 

(2) (Quasi) R-matrix \cite{Dr87, Lus93}

(3) Canonical basis  \cite{Lus90, Lus92, Lus93} \cite{Ka91}

(4) Quantum Schur duality \cite{Jim86}

(5) Super type A Kazhdan-Lusztig theory \cite{Bru03} \cite{CLW15} \cite{BLW17}

(6) Hall algebra \cite{Rin90} \cite{Br13} 

(7)
Current presentation of affine quantum groups \cite{Dr87} \cite{Be94, Da12}

(8)
Braid group action \cite{Lus90, Lus93}

(9)
Categorification \cite{KLI, R08}


(+) $\circ$ $\circ$ $\circ$ $\circ$ $\circ$ $\circ$

We apologize beforehand for omitting many important constructions for quantum groups in the above list, and your favorite construction may likely fall into the black holes in Item (+). Also, the references listed above are mostly samples of the original contributions, and there are often dozens or hundreds of additional works which are not cited. 

 The list (1)--(9) are so arranged  that they are to be matched with the  $\imath$-generalizations which we shall describe later in the same ordering.
\subsection{Quantum symmetric pairs}

Recall a symmetric pair $(\g, \g^\theta)$ consists of a semisimple Lie algebra $\g$ and an involution $\theta$ on $\g$. The classification of irreducible symmetric pairs is equivalent to the classification of real forms of complex simple Lie algebras (which goes back to \'E.~ Cartan); for example, they can now be classified in terms of Satake diagrams. 

Let $\bU =\bU_v(\g)$ be a quantum group (of finite type) with comultiplication $\Delta$ in  \eqref{eq:Delta}. 
According to Gail Letzter \cite{Let99, Let03, Let19}, a quantum symmetric pair $(\bU, \Ui)$ consists of a Drinfeld-Jimbo quantum group $\bU$ and its (right) coideal subalgebra $\Ui$ (i.e., $\Delta: \Ui \rightarrow \Ui \otimes \bU$) which specializes at $v \mapsto 1$ to $U(\g^\theta)$. Starting with the Satake diagrams, Letzter constructed the corresponding quantum symmetric pairs. The $\Ui$ comes with parameters, which reflects the fact that there is a family of (explicit) embeddings of $\Ui$ into $\bU$; see \eqref{eq:embed}. A generalization of quantum symmetric pairs of Kac-Moody type was carried out by Kolb in \cite{Ko14}, now a standard reference in the subject; Kolb's conventions are compatible with those in Lusztig's book \cite{Lus93}. As we may deal with $\Ui$ alone, we shall call $\Ui$ an {\em $\imath$quantum group}.

Letzter's foundational work on quantum symmetric pairs was motivated by harmonic analysis on quantum symmetric spaces, generalizing earlier examples given by Koornwinder, Gavrilik-Klimyk, Noumi and others. Letzter's work was ahead of her time.
\subsection{Goal}

We take the view that $\imath$quantum groups are a vast generalization of quantum groups (as real reductive groups are a generalization of compact or complex reductive groups). 
\begin{example}
[Quantum groups as $\imath$quantum groups]
\label{ex:diag}
Consider the diagonal symmetric pair $(\g \times \g, \g^\triangle)$, where $\g^\triangle$ is a diagonal copy of $\g$. Similarly, we have a quantum symmetic pair of diagonal type $(\bU \otimes \bU, \bU)$; the embedding is given by 
$(\omega \otimes 1) \Delta: \bU \rightarrow \bU \otimes \bU$ (see \cite{BW18b}), where one checks the image of $\bU$ is a coideal subalgebra of $\bU \otimes \bU$.
\end{example}

The goal of this report is to survey some recent advances on quantum symmetric pairs and $\imath$quantum groups, generalizing Items (1)-(9) above. The $\imath$-analogs of the constructions in Items (2)-(5) were initiated by Huanchen Bao and the author in \cite{BW18a}. It was proposed loc. cit. that all fundamental (algebraic, geometric, categorical) constructions in quantum groups should be generalized to $\imath$quantum groups. 

The good news is that all Items (1)-(9) admit genuine $\imath$-generalizations, while the bad news, or the exciting news for an optimist, is that many $\imath$-generalizations are not yet in full generality. 
A reader might well be tempted to try one's hands in developing these $\imath$-generalizations in greater generality.
She is encouraged to pick her own favorite construction in quantum groups in Item~(+); even better, to supply its missing $\imath$-generalization. 

All the known constructions indicate that the $\imath$-generalizations (when done right!) look natural and inevitable, though proofs are often much more lengthy and challenging. The complications in $\imath$quantum groups are often caused by 

(i) absence of triangular decompositions; 

(ii) presence of many rank 1 types \emph{and} parameters;

(iii) presence of nonuniform Serre type relations;

(iv) hidden (non-obvious) integral forms.

\subsection{A quick overview}

A Serre type presentation for an $\imath$quantum group $\Ui$ was due to Letzter \cite{Let03} in finite type, and has been generalized since then to Kac-Moody type in various forms (cf. \cite{Ko14} \cite{BK15, CLW18, dC19, KY21}). This can be viewed as an $\imath$-analog of the construction in Item (1).  

Let us comment in some details on the constructions in Items (2)-(5). 

As a generalization of Jimbo-Schur duality \cite{Jim86}, an $\imath$Schur duality between a  {\em quasi-split} type AIII $\imath$quantum group $\Ui$ and type B Hecke algebra was formulated in \cite{BW18a} (and \cite{Bao17, BWW18}). Recently, the Jimbo-Schur duality and quasi-split $\imath$Schur duality has been uniformly generalized using a general $\imath$quantum group of type AIII by Yaolong Shen and the author in \cite{SW21}. This has led unexpectedly to quasi-parabolic Kazhdan-Luszig bases associated to (possibly non-parabolic) reflection subgroups of the type B Weyl group, extending the classic constructions of Kazhdan-Lusztig \cite{KL79} and Deodhar \cite{De87}.

The $\imath$canonical basis and $\imath$Schur duality (in quasi-split type AIII) were motivated by and played a key role in formulating a super type B Kazhdan-Lusztig theory (which was a decades old open problem) \cite{BW18a}; see also Bao \cite{Bao17} for a super type D formulation. Our approach was inspired by Brundan-Kazhdan-Lusztig conjecture in super type A \cite{Bru03} and its proof by Cheng, Lam and the author \cite{CLW15}. 

The $\imath$-analog of quasi R-matrix, known as quasi K-matrix nowadays, was formulated by Bao and the author in \cite{BW18a} as a key step in the construction of $\imath$canonical bases, and a proof for its existence in great generality has been given by Balagovic-Kolb in \cite{BK19}; see Appel-Vlaar \cite{AV20} for a more recent reformulation and generalization. The quasi K-matrix further leads to a construction of the K-matrix \cite{BW18a} \cite{BK19}, which was shown in \cite{BK19} to provide solutions to the reflection equation.

As an extension of Lusztig's approach \cite{Lus90, Lus92}, a theory of $\imath$canonical bases arising from quantum symmetric pairs has been systematically developed by Bao and the author in \cite{BW18b, BW21}. The $\imath$canonical bases on based $\bU$-modules (viewed as $\Ui$-modules) and on the modified $\imath$quantum groups are established; also see \cite{BKLW, LiW18} and \cite{FL+20} for a geometric approach in type AIII. The $\imath$canonical basis in split rank 1, also known as $\imath$divided powers \cite{BW18a, BeW18}, has found applications in works with Xinhong Chen and Ming Lu  \cite{CLW18, CLW21, LW20a}. H.~Watanabe \cite{Wat21a, Wat21} has  developed a crystal approach (\`a la Kashiwara \cite{Ka91}) to $\imath$canonical bases of $\Ui$-modules, for some quasi-split finite types. 

Motivated by earlier constructions of Bridgeland \cite{Br13} and then M.~Gorsky \cite{Gor18} (extending the foundational work of Ringel \cite{Rin90}), Lu-Peng and Lu developed semi-derived Hall algebras associated to 1-Gorenstein algebras (see \cite{LP16} and \cite[Appendix A]{LW19a}). Lu and the author introduced in \cite{LW19a, LW20a} the $\imath$Hall algebras, i.e., semi-derived Hall algebras associated to the $\imath$quiver algebras, to realize the (\emph{universal}) $\imath$quantum groups $\tUi$. This is a conceptual $\imath$-analog of (6), generalizing \cite{Br13}. 


In \cite{LW20b}, Lu and the author has formulated a Drinfeld type presentation for affine $\imath$quantum groups of split ADE type; this has been generalized to split BCFG type by Weinan Zhang \cite{Zh21}. This provides an $\imath$-analog of (7). 
 Lu and Ruan are developing $\imath$Hall algebras of the weighted projective lines \cite{LR21} to realize the affine $\imath$quantum groups in the new presentation, generalizing the rank one construction in \cite{LRW20}.  

The braid group actions associated to relative Weyl groups on (mostly quasi-split) finite type $\Ui$ were obtained by Kolb-Pellegrini in \cite{KP11} (via computer computation), where the existence of such an action on an arbitrary $\imath$quantum group of finite type was conjectured. Relative braid group actions on $\tUi$ (of quasi-split Kac-Moody type) were obtained in \cite{LW21a} via reflection functors in $\imath$Hall algebras, and it becomes clear that the universal $\imath$quantum groups provide a conceptual framework for braid group actions. We announce an intrinsic approach developed by W.~ Zhang and the author in \cite{WZ21} to relative braid group action on $\tUi$ of {\em arbitrary} finite type (and on $\bU$-modules) with explicit formulas, providing a conceptual $\imath$-analog of (8). The quasi K-matrix plays a basic role in our formulation.  

This survey is organized as follows. 
The $\imath$-counterparts of the quantum group highlights (1)-(8) will be formulated in Sections \ref{sec:def}--\ref{sec:affine} below.
In Section~\ref{sec:open}, we discuss additional topics in $\imath$quantum groups including the $\imath$-analog of (9), and list some open problems.

\section{Quantum symmetric pairs: definition}
  \label{sec:def}

\subsection{Quantum groups} 

We mostly follow notations in the book \cite{Lus93}. Denote by $\N$ the set of non-negative integers. Let $(Y,X, \langle\cdot,\cdot\rangle, \cdots)$ be a root datum of type $(\I, \cdot)$; cf. \cite{Lus93}. The quantum group $\bU$ associated with this root datum  $(Y,X, \langle\cdot,\cdot\rangle, \cdots)$ is the associative $\Qq$-algebra generated by $E_{i}$, $F_{i}$ for $i \in \I$ and $K_{\mu}$ for $\mu \in Y$, subject to standard relations which can be found in \cite{Lus93}. 
 %
 %
%
Let $W$ denote the Weyl group generated by simple reflections $s_i$ for $i \in \I$. The $\Qq$-algebra $\bU$ admits a Chevalley involution $\omega$ such that
\begin{align}
\label{inv:omega}
\omega(E_i) =F_i, \quad
\omega(F_i) =E_i, \quad
\omega(K_\mu) =K_{-\mu},
\quad
\text{ for }i\in \I, \mu \in Y.
\end{align}
For any $i \in \I$, we set $v_{i} = v^{\frac{i \cdot i}{2}}$. For $i\in \I, n, s \in \N$ with $0 \le s \le n$, we define 
$[n]_i = \frac{v_i^n-v_i^{-n}}{v_i - v^{-1}_i}$ and $[s]^!_i = \prod^s_{k=1} [k]_i$, and  
$\qbinom{n}{s}_i  = \frac{[n]^!_i}{[s]^!_i [n-s]^!_i}.$ 

Let $\bU^+$, $\bU^0$ and $\bU^-$ be the $\Qq$-subalgebra of $\bU$ generated by $E_{i} (i \in \I)$, $K_{\mu} (\mu \in Y)$, and $F_{i}  (i \in \I)$, respectively. 
Let $\mA =\Z[v,v^{-1}].$ We let $_\mA\bU^-$ (respectively, $_\mA\bU^+$) denote the $\mA$-subalgebra of $\bU^-$ (respectively, $\bU^+$) generated by all divided powers $F^{(a)}_{i}=F_{i}^s/[s]_i^!$ (respectively, $E^{(a)}_{i}$). With $\Ktilde_{\pm i} := K_{\pm \frac{i \cdot i}{2} i}$, the coproduct $\Delta: \bU \rightarrow \bU \otimes \bU$ is given by
\begin{equation}\label{eq:Delta}
\Delta(E_i)  = E_i \otimes 1 + \tilde{K}_i \otimes E_i, \quad \Delta(F_i) = 1 \otimes F_i + F_i \otimes \tilde{K}_{-i}, \quad \Delta(K_{\mu}) = K_{\mu} \otimes K_{\mu}.
\end{equation}

%

Let 
$X^+ = \big\{\lambda \in X \mid \langle i, \lambda \rangle \in {\N}, \forall i \in \I \big \}$ 
 be the set of dominant integral weights. 
 By $\la \gg 0$ we shall mean that the integers $\langle i, \lambda \rangle$ for all $i$ are sufficiently large. The Verma module $M(\lambda)$ of highest weight $\la \in X$ has a unique simple quotient $\bU$-module $L(\lambda)$ with a highest weight vector $\eta_\lambda$. We define a $\bU$-module $^\omega L(\lambda)$, which has the same underlying vector space as $L(\lambda)$ but with the action twisted by the Chevalley involution $\omega$ in \eqref{inv:omega}; then $^\omega L(\lambda)$ is simple of lowest weight $-\la$ with lowest weight vector denoted by $\xi_{-\la}$.  
 For $\la \in X^+$, we let $_\mA L(\lambda) ={_\mA\bU^-}  \eta_\la$ and $^\omega _\mA L(\lambda) ={_\mA\bU^+ } \xi_{-\la}$ be the $\mA$-submodules of $L(\lambda)$ and $^\omega L(\lambda) $. 

%
There is a canonical basis $\B$ on $\bU^-$, inducing a canonical basis on $\bU^+$ via the standard isomorphism $\bU^-\cong \bU^+$. For each $\la\in X^+$, there is a subset $\B(\la)$ of $\B$ so that $\{b^- \eta_\la \mid b\in \B(\la)\}$  forms a canonical basis of $L(\la)$. For $w\in W$, let $\eta_{w\la}$ denote the unique canonical basis element of weight $w\la$. 

Let $\Udot =\oplus_{\zeta \in X} \Udot \one_\zeta$ be the modified quantum group and $\aA \Udot$ be its $\A$-form. Then $\Udot$ admits a canonical  basis 
$\dot{\B}  = \{ b_1 \diamondsuit_{\zeta} b_2   \vert (b_1, b_2) \in \B \times \B, \zeta \in X \}$, compatible with canonical bases on $^\omega L(\lambda) \otimes L(\mu)$, for $\lambda, \mu \in X^+$; cf. \cite[Part ~IV]{Lus93}.
For any $\Iblack \subset \I$, let $\bU_{\Iblack}$ be the $\Qq$-subalgebra of $\bU$ generated by $F_{i}$, $E_{i}$ and $K_{i} ( i \in \Iblack)$.  
Let $\B_{\I_\bullet}$ be the canonical basis of $\bU_{\Iblack}^-$. 

%
%
\subsection{Satake diagrams and admissible pairs}
  \label{subsec:adm}
  
 Let $\tau$ be an involution of the Cartan datum $(\I, \cdot)$; we allow $\tau =\Id$. We further assume that $\tau$ extends to  an involution on $X$ and an involution on $Y$, respectively, such that the perfect bilinear pairing is invariant under the involution $\tau$. Given a  finite type $\Iblack \subset \I$,  let $W_{\I_\bullet} =\langle s_i \mid i\in \Iblack\rangle$ be the parabolic subgroup of $W$ with longest element $w_{\bullet}$, and let $\rho_{\bullet}$ (and $\rho^{\vee}_{\bullet}$) be the half sum of all positive roots (and coroots) in the root system $R_{\bullet}$ (and $R^{\vee}_{\bullet}$). Set $\I_{\circ} = \I \backslash \I_{\bullet}.$
A pair $(\I =\I_{\bullet} \cup \Iwhite, \tau)$ is called {\em admissible} (cf. \cite[Definition~2.3]{Ko14}) if 	$\tau (\I_{\bullet}) = \I_{\bullet}$,
the actions of $\tau$ and $-w_{\bullet}$ on $\I_{\bullet}$ coincide, and 
$\langle \rho^\vee_{\bullet}, j' \rangle \in \Z$ whenever $\tau j =j \in \I_{\circ}$. We regard admissible pairs as a synonyms for Satake diagrams. 


Note that $\inv = -w_{\bullet} \circ \tau$ is an involution of $X$ and $Y$. For any $\lambda \in X$ (or $Y$), we shall write $\lambda^\tau = \tau(\lambda), \lambda^\inv = \inv(\lambda)$.
 Following \cite{BW18b}, we introduce the $\imath$-weight and $\imath$-root lattices
\begin{align*}
X_{{\imath}} &= X /  \breve{X},   \quad \text{ where } \; \breve{X}  = \{ \la - \la^\inv \vert \la \in X\},
 \\
Y^{\imath} &= \{\mu \in Y \big \vert \mu^\inv =\mu \}.
\end{align*}
For $\la \in X$, we denote its image in $X_{\imath}$ by $\overline{\la}$.
The involution $\tau$ of $\I$ induces an isomorphism of the $\Qq$-algebra $\bU$, denoted also by $\tau$,
which sends $E_i \mapsto E_{\tau i}, F_i \mapsto F_{\tau i}$, and $K_\mu \mapsto K_{\tau \mu}$. 
\subsection{Quantum symmetric pairs}

We review the definition of quantum symmetric pair $(\bU, \Ui)$, where $\Ui$ is a coideal subalgebra of $\bU$, following \cite{Ko14}; also see \cite{BK15, BK19, BW18b}. Letzter's convention is a little different. 

An {\em $\imath$quantum group} $\Ui$ is the $\Qq$-subalgebra of $\bU$ generated by 
\begin{align}
B_i: = F_{i}  &+ \vs_i T_{w_{\bullet}} (E_{\tau i}) \tilde{K}^{-1}_i 
+ \kappa_i \tilde{K}^{-1}_{i} \; (i \in \I_{\circ}), 
 \label{eq:Bi}  
\quad K_{\mu} \,(\mu \in Y^{\imath}), \quad F_i, E_{i} \,(i \in \I_{\bullet}),
\end{align}
where $
\vs_{i}  \in \Qq^\times,   \kappa_i \in \Qq,  \text{ for }  i \in \I_{\circ},
$ 
are parameters, and $T_w= T''_{w,+1}$ denotes a braid group operator as in \cite{Lus93}. 
Denote by $\bU^{\imath 0}$ the $\Qq$-subalgebra of $\Ui$ generated by $K_{\mu}$ $(\mu \in Y^{\imath})$,  
and denote the embedding via \eqref{eq:Bi} by
\begin{align}
 \label{eq:embed}
 \imath: \Ui \longrightarrow \bU, \quad x \mapsto x^\imath.
\end{align}
The parameters $\vs_i, \kappa_i$ are required to satisfy Conditions \eqref{kappa}-\eqref{vs=}: 
\begin{align}
 \label{kappa}
 \begin{split}
 \kappa_i & 
 =0 \; \text{ unless } ``\tau i =i, \langle i, j' \rangle = 0 \; \forall j \in \Iblack,
\&  \langle k,i' \rangle \in 2\Z \; \forall \tau k =k =w_\bullet k \in \Iwhite";
\end{split}
\\
\vs_{i} & =\vs_{{\tau i}} \text{ if }    i \cdot \theta (i) =0.
\label{vs=}
\end{align}
The conditions \eqref{kappa}--\eqref{vs=} ensure the quantum Iwasawa decomposition and hence $\Ui$ has the expected size \cite{Let03} \cite{Ko14}. 
By definition, the algebra $\Ui$ contains $\bU_{\Iblack}$ as a subalgebra.

In the remainder of this paper, we shall impose the following additional conditions on parameters besides \eqref{kappa}--\eqref{vs=}, as required for the construction of quasi K-matrix and $\imath$canonical basis.
\begin{align}
  \label{parameters}
\vs_{i} &,  \kappa_i \in \Z[v,v^{-1}], \qquad \text{ for }  i \in \I_{\circ},
\\
\overline{\kappa_i} &= \kappa_i,   \qquad 
\vs_{{\tau i}} = (-1)^{ \langle 2\rho^\vee_{\bullet},  i' \rangle } v_i^{-\langle i, 2\rho_{\bullet}+\wb\tau i ' \rangle} {\ov{\vs_{i} }}.   \label{vs2}
\end{align}
Conditions~\eqref{vs2} appeared in \cite{BK19}.

The nontrivial Serre presentation of $\Ui$ was given by Letzter \cite{Let03} in finite type, and further generalized in different forms \cite{Ko14} \cite{BK15, CLW18, KY21}. 
Certain Serre-Lusztig (or higher order Serre) relations for $\Ui$ have been obtained in \cite{CLW21} via $\imath$divided powers (see Example~\ref{ex:rank1}). 

\begin{example} 
We call $\Ui$ \emph{quasi-split} if $\Iblack =\emptyset$ and \emph{split} if in addition $\tau =\Id$. A split $\imath$quantum group $\Ui$ is the subalgebra of $\bU$ generated by $B_i =F_{i}  + \vs_i E_{i} \tilde{K}^{-1}_i + \kappa_i \tilde{K}^{-1}_{i}$ $(i\in \I)$; often $\kappa_i=0$ thanks to \eqref{kappa}. Relations in split $\Ui$ are given in \eqref{eq:S1}--\eqref{eq:S2} (with $\K_i= -v_i^2 \vs_i$) for ADE type and in \eqref{iDRG5'}--\eqref{iDRG6'} (with $k=l=0$) for BCFG type. 
\end{example}

\section{(Quasi) K-matrices}

As predicted in \cite{BW18a} and established in \cite{BK15, Ko21} (also cf. \cite{CLW18}), 
there is a unique anti-linear bar involution of the $\Q$-algebra $\Ui$, denoted by $\psi_{\imath}$, such that $\psi_{\imath} (v) =v^{-1}$ and 
\begin{align*}
 \psi_{\imath} (B_i) = B_i \; (i \in \Iwhite), \quad \ipsi(F_j) = F_j, \quad \ipsi(E_j) = E_j  \;(i \in \I_{\bullet}), \quad \ipsi (K_\mu) = K_{-\mu} \; (\mu \in Y^{\imath}).
\end{align*}

The formulation of the {\em quasi $K$-matrix} $\Upsilon$ below (called an {\em intertwiner} earlier) was due to Bao and the author \cite{BW18a}; its existence in great generality has been established in \cite{BK19} (also cf. \cite[Remark 4.9]{BW18b}) with additional technicality removed in \cite{BW21, Ko21}. The quasi $K$-matrix for quantum symmetric pair of diagonal type reduces to the quasi $R$-matrix (cf. \cite{Lus93}).

\begin{theorem}  \cite{BW18a} \cite{BK19} \cite{BW18b}
   \label{thm:Upsilon}
There exists a unique family of elements $\Upsilon_{\mu} \in \bU^+_{\mu}$, for $\mu \in \N \I$, 
such that $\Upsilon_0 =1$ and $\Upsilon = \sum_{\mu \in \N \I} \Upsilon_{\mu}$ satisfies the following identity:
\begin{equation}\label{eq:Upsilon}
 \ipsi(u) \Upsilon = \Upsilon \psi(u), \qquad \text{for all } u \in \Ui.
\end{equation}
Moreover, $\Upsilon_{\mu} =0$ unless ${\mu^\inv} = - \mu \in X$. (Recall $\inv = -w_{\bullet} \circ \tau$.)
\end{theorem}

Quasi K-matrix was originally formulated in order to construct a new bar involution and then $\imath$canonical bases on based $\bU$-modules (see \S\ref{sec:iCB} below). On the other hand, a suitable twisting of $\Upsilon$ by elements in $\bU^{\imath 0}$ provides certain $\Ui$-module isomorphisms \cite[\S2.5]{BW18a} \cite{BK19} (also see \cite[\S4.5]{BW18b} for a different formulation), nowadays known as the K-matrix. It is shown in \cite{BK19} that the K-matrix provides solutions to the reflection equation, an $\imath$-analog of the Yang-Baxter equation. 
There has been a reformulation of quasi K-matrix in \cite{AV20} without referring explicitly to the bar involution of $\Ui$; this has the advantage of making sense of quasi K-matrices in general parameters satisfying \eqref{kappa}--\eqref{vs=}. 

\section{Canonical bases arising from quantum symmetric pairs}
 \label{sec:iCB}

\subsection{Based modules}

The quasi $K$-matrix $\Upsilon$ in (a completion of) $\bU^+$ induces a $\Qq$-linear map 
$
\Upsilon: M \otimes L(\la) \rightarrow M \otimes L(\la),
$ 
for any $\la \in X^+$ and any weight $\bU$-module $M$ with weights bounded above. A $\Ui$-module $M$ equipped with an anti-linear involution $\ipsi$ is called {\em $\imath$-involutive} if $\ipsi(u m) = \ipsi(u) \ipsi(m),$ for all $u \in \Ui, m \in M.$
Let $(M,B)$ be a based $\bU$-module \cite[IV]{Lus93} with weights bounded above. We denote the bar involution on $M$ by $\psi$. Then $M$ is an $\imath$-involutive  
$\Ui$-module with involution (see \cite{BW18a})
\begin{equation*}
\ipsi := \Upsilon \circ \psi.
\end{equation*}

Assume that $(M,B)$ is a based $\bU$-module with weights bounded above such that $\Upsilon$ preserves the $\A$-submodule $\aM$. Then the  $\Qq$-linear map $\ipsi = \Upsilon \circ \psi$ (and subsequently, $\Upsilon$) preserves the $\A$-submodule $\aM \otimes_\A \aL(\la)$, for any $\la\in X^+$; see \cite{BW21}. 
In particular, the involution $\ipsi$ on the $\imath$-involutive $\Ui$-module $L(\la_{1})\otimes \ldots \otimes L(\la_\ell)$ preserves the $\mA$-submodule ${}_\mA L(\la_{1})\otimes_{\mA} \ldots \otimes_{\mA} {}_\mA L(\la_\ell)$, where  $\lambda_i \in X^+$ for $1 \le i \le \ell$.
For $(\bU, \Ui)$ of finite type, a stronger statement holds \cite{BW21}: $\Upsilon_\mu \in \aA \bU^+$, for each $\mu$. This generalizes the integrality of quasi R-matrix of finite type due to Lusztig \cite[24.1.6]{Lus93}. 

Define a partial order $\leq$ on the weight lattice $X$ such that $\lambda \le \lambda' \text{ if and only if } \lambda' -\lambda \in \N \I.$ For an element $x$ in $\bU$ or in a $\bU$-module $M$ of weight $\mu \in X$, we write $|x| =\mu$.

\begin{theorem}  \cite{BW18b, BW21}
\label{thm:iCBmodule}
Let $(M,B)$ be a based $\bU$-module with weights bounded above. Assume the involution $\ipsi =\Upsilon \psi$ of $M$ preserves the $\A$-submodule $\aM$. Then, 
\begin{enumerate}
\item
the $\Ui$-module $M$ admits a unique basis (called $\imath$canonical basis)
$
B^\imath  := \{b^\imath \mid b \in B \}
$
which is $\ipsi$-invariant and of the form
$
b^\imath \in b +\sum_{b' \in B, |b| <  |b'|}
v^{-1}\Z[v^{-1}] b'.
$ 

\item
$B^\imath$ forms an $\mA$-basis for the $\mA$-lattice ${}_\mA M$ (generated by $B$), and
forms a $\Z[v^{-1}]$-basis for the $\Z[v^{-1}]$-lattice $\mc{M}$ (generated by $B$).
\end{enumerate}
\end{theorem}  
In particular, $L(\la_{1})\otimes \ldots \otimes L(\la_\ell)$, where  $\lambda_i \in X^+$ for all $i$, admits an $\imath$canonical basis.
Theorem~\ref{thm:iCBmodule} can be further generalized in \cite{BWW20, BW21} to provide an $\imath$canonical basis on $N\otimes L(\la)$, for a based $\Ui$-module $N$ and $\la \in X^+$.

\subsection{Canonical bases on modified $\imath$quantum groups}

Following \cite[IV]{Lus93}, we can define a modified $\imath$quantum group $\Uidot$ (an associative $\Q(q)$-algebra structure without unit) such that 
$\Uidot = \bigoplus_{\lambda \in X_\imath} \Ui \mathbf{1}_{\lambda}$; see \cite{BW21}. In contrast to quantum groups, the $\mA$-form of $\Uidot$ is far from being obvious (even in rank 1).

For $\lambda, \mu \in X^+$ and $w\in W$, denoting $\eta_\lambda^\bullet :=\eta_{w_\bullet \lambda}$, we introduce the following $\bU$-submodule:
\begin{align*}
L^\imath(\la,\mu) := \bU (\eta_{w_\bullet \lambda} \otimes \eta_{\mu}) \subset L(\lambda) \otimes L(\mu), 
\end{align*}
which can be shown to be a based $\bU$-module and hence admits an $\imath$canonical basis by Theorem~\ref{thm:iCBmodule}. 
Let $\zeta = \wb \lambda +\mu$ and $\zeta_\imath = \overline{\zeta}$. The following hold \cite{BW18b, BW21}:

$\triangleright$		The $\imath$canonical basis of $L^{\imath}(\lambda, \mu)$ is of the form $\B^\imath(\la, \mu) =\big\{ (b_1 \diamondsuit_{\zeta_\imath} b_2 )_{w_\bullet\la,\mu}^{\imath} \vert 
	(b_1, b_2) \in \B_{\Iblack} \times \B \big\} \backslash \{0\}$, where 
 $(b_1 \diamondsuit_{\zeta_\imath} b_2 )_{w_\bullet\la,\mu}^{\imath}$ is $\ipsi$-invariant and lies in 
	\begin{equation*}  
	 (b_1 \diamondsuit_{\zeta} b_2 ) (\eta_\lambda^\bullet \otimes \eta_{\mu})    +\!\!\! \sum_{|b_1'|+|b_2'| \le |b_1|+|b_2|} \!\!\!v^{-1}\Z[v^{-1}] (b'_1 \diamondsuit_{\zeta} b'_2 ) (\eta_\lambda^\bullet \otimes \eta_{\mu}). 
	\end{equation*}

$\triangleright$	We have a projective system $\big \{ L^\imath(\lambda+\nu^\tau, \mu+\nu)  \big \}_{\nu \in X^+}$ of $\Ui$-modules, where 
\begin{equation*} 
\pi_{\nu+ \nu_1, \nu_1}: L^{\imath} (\lambda+\nu^\tau+\nu_1^\tau, \mu+\nu+\nu_1) \longrightarrow L^{\imath} (\lambda+\nu^\tau, \mu+\nu), \quad \nu, \nu_1 \in X^+,
\end{equation*}
is the unique homomorphism of $\Ui$-modules mapping $\etab_{\lambda+\nu^\tau+\nu_1^\tau} \otimes \eta_{\mu+\nu+\nu_1}$ to $\etab_{\lambda+\nu^\tau} \otimes \eta_{\mu+\nu}.$ The $K$-matrix \cite{BW18a} \cite{BK19} \cite{BW18b} plays a role here. 


The following theorem generalizes \cite[Chap. 25]{Lus93}. For quantum symmetric pairs of diagonal type in Example~\ref{ex:diag} or a trivial pair $(\bU,\bU)$, it reduces to Lusztig's setting.

\begin{theorem}  \cite{BW18b}  \cite[Theorem 7.2]{BW21}
  \label{thm:iCBUi}
		Let ${\zeta_\imath} \in X_{\imath}$ and $(b_1, b_2) \in B_{\Iblack} \times B$. 
	\begin{enumerate}
		\item	
		There is a unique element $u =b_1 \diamondsuit^\imath_{\zeta_\imath} b_2  \in \Uidot$ such that 
			\[
				u (\eta^{\bullet}_\lambda \otimes \eta_\mu) 
				= (b_1 \diamondsuit_{\zeta_\imath} b_2)_{\wb\lambda,\mu}^\imath				\in L^{\imath} (\lambda, \mu), 
\text{ for all } \lambda, \mu \gg0 \text{ with }\overline{\wb \lambda+\mu} = {\zeta_\imath}.
\]
		\item	The element $b_1 \diamondsuit^\imath_{\zeta_\imath} b_2$ is $\ipsi$-invariant.
		\item	The set $\Bdot^\imath =\{ b_1 \diamondsuit^\imath_{\zeta_\imath} b_2 \mid {\zeta_\imath} \in X_\imath, (b_1, b_2) \in B_{\Iblack} \times B \}$ 
		forms a $\Qq$-basis of $\Uidot$.
	\end{enumerate}
\end{theorem}

\begin{example} [$\imath$Divided powers] 
 \label{ex:rank1}
Consider the quantum symmetric pair of split rank one $(\bU, \Ui) =(\bU_v(\mathfrak{sl}_2), \Qq [B_i])$ associated to $\I=\{i\}$, via the embedding $\Ui \rightarrow \bU$, $B_i \mapsto F_i + \vs_i E_i K_i^{-1}$. It is a new phenomenon of $\imath$quantum groups \cite{BW18a} that there are two different modified forms of $\Ui$, denoted by $\dot \bU^{\imath} {\mathbf 1}_{\bar 0}$ and $\dot \bU^{\imath} {\mathbf 1}_{\bar 1}$, depending on a parity $X_\imath = \{\bar 0, \bar 1\}$, which are compatible with the parity of highest weights of finite-dimensional simple $\bU$-modules.

We define the {\em $\imath$divided powers} of $B$ to be
\begin{eqnarray*}
&&
B_{i,\bar 1}^{(m)}=\frac{1}{[m]_{i}^!}\left\{ 
\begin{array}{ccccc} B_i\prod_{j=1}^\ell (B_i^2-v_i\vs_i[2j-1]_{i}^2 ) & \text{if }m=2\ell+1,\\
\prod_{j=1}^\ell (B_i^2-v_i\vs_i[2j-1]_{i}^2) &\text{if }m=2\ell; \end{array}\right.
  \\
&&
B_{i,\bar 0}^{(m)}= \frac{1}{[m]_{i}^!}\left\{ \begin{array}{ccccc} B_i\prod_{j=1}^\ell (B_i^2-v_i\vs_i[2j]_{i}^2 ) & \text{if }m=2\ell+1,\\
B_i^2 \prod_{j=1}^{\ell-1} (B_i^2-v_i\vs_i[2j]_{i}^2) &\text{if }m=2\ell. \end{array}\right.
\end{eqnarray*}
These formulas (with parity swapped) appeared first in \cite[Conjecture~ 4.13]{BW18a} (in terms of $B_i$ in \eqref{eq:Bi} with $\kappa_i=1, \vs_i=v_i^{-1}$); see \cite{CLW18} for application to Serre relations for $\imath$quantum groups.

Set the parameter $\vs_i=v_i^{-1}$. Then ${\mathbf B}_{\bar 0} :=\{B_{i,\bar 0}^{(m)} \mid m \ge 0 \}$ (and respectively, ${\mathbf B}_{\bar 1} :=\{B_{i,\bar 1}^{(m)} \mid m \ge 0 \}$) forms the $\imath$canoical basis for the modified $\imath$quantum group $\dot \bU^{\imath} {\mathbf 1}_{\bar 0}$ (and respectively, $\dot \bU^{\imath} {\mathbf 1}_{\bar 1}$). 
Let $L(n)$ be the simple $\bU$-module of highest weight $n\in \N$ with highest weight vector $\eta$. Then, for $n$ even,  ${\mathbf B}_{\bar 0}\eta =\{B_{i,\bar 0}^{(m)} \eta \mid n\ge m \ge 0 \}$ forms the $\imath$canonical basis for $L(n)$ and $B_{i,\bar 0}^{(m)} \eta=0$, for $m>n$; similar claims hold for $L(n)$ with $n$ odd and ${\mathbf B}_{\bar 1}$. See \cite{BeW18}. 
\end{example}

\section{Quantum Schur dualities}

\subsection{Quasi-parabolic Kazhdan-Lusztig bases}

Let $W_d$ be the Weyl group of type $B_d$ generated by simple reflections $s_i$, for $0\le i \le d-1$. 
It contains the symmetric group $\mathfrak S_d$ as a subgroup generated by $s_i$, for $1\le i \le d-1$. 
For $p \in v^\Z$, the Hecke algebra $\HB$ of type $B_d$ is a $\Qq$-algebra generated by $H_i$ $(0 \le i \le d-1)$ such that $(H_0+p^{-1})(H_0-p)=0$, $(H_i+v^{-1})(H_i-v)=0$ for $i\ge 1$, and braid relations hold; $\HB$ admits a bar involution $\bar{\phantom{x}}$ such that $\ov{v} = v^{-1}$ and $\ov{H_i}= H_i^{-1}$, for all $i$. 

For $x\in \mathbb R$ and $m \in \N$, we denote $[x.. x+m] =\{x, x+1, \ldots, x+m \}$. 
For $a \in \Z_{\ge 1}$, we denote by $
\II_a = \big[\frac{1-a}2 .. \frac{a-1}2 \big].$ 
For $r, m \in \N$ (not both zero), we introduce a new notation $\Ibw : =\II_{2r+m}$ to indicate a fixed set partition:
$\Ibw =\Iwl \cup \Ib \cup \Iwr$, where 
\begin{align}
 \label{eq:III}
\Iwr = \Big [\frac{m+1}{2} .. r+\frac{m-1}{2} \Big ], 
\qquad
\Ib = \Big [\frac{1-m}{2} .. \frac{m-1}{2} \Big ], 
\qquad
\Iwl = - \Iwr.
\end{align}
We view $f  \in \Ibw^d$ as a map $f: \{1, \ldots, d\} \rightarrow \Ibw$, and identify 
$f=(f(1), \ldots, f(d))$. 
We define a right action of $W_d$ on $\Ibw^d$ such that,
for $f\in \Ibw^d$ and $0\leq j\leq d-1$,
\begin{equation*}
f\cdot s_j = 
\begin{cases} 
 (\cdots,f(j+1),f(j),\cdots),&\text{ if } j>0, \\
 (-f(1),f(2),\cdots,f(d)), &\text{ if } j=0,\ f(1)\in \Iwl\cup \Iwr, \\
 (f(1),f(2),\cdots,f(d)), &\text{ if } j=0,\ f(1)\in \Ib.
\end{cases} 
\end{equation*} 

Let $p \in v^\Z$. Consider the $\Qq$-vector space $\V=\bigoplus_{a\in \Ibw}\Q(v) u_a.$
Given $f=(f(1), \ldots, f(d)) \in \Ibw^d$, we denote
$M_f = u_{f(1)} \otimes u_{f(2)} \otimes \ldots \otimes u_{f(d)}.$
We shall call $f$ a weight and $\{M_f \mid f\in \Ibw^d \}$ the standard basis for $\V^{\otimes d}$. 
There is a right action of the Hecke algebra $\Hy_{B_d}$ on $\V^{\otimes d}$ as follows (see \cite{SW21}):
\[
M_f\cdot H_i=\left\{
\begin{aligned}
&M_{f\cdot s_i}+(v-v^{-1})M_{f},\ \ & \text{ if } f(i)<f(i+1),\ i>0, \\
&M_{f\cdot s_i},\ \ &\text{ if } f(i)>f(i+1),\ i>0, \\
&vM_f,\ \ & \text{ if } f(i)=f(i+1),\ i>0, \\
&M_{f\cdot s_i}+(p-p^{-1})M_f, \ \ &\text{ if } f(1)\in \Iwr,\ i=0, \\
&M_{f\cdot s_i},\ \ &\text{ if } f(1)\in \Iwl,\ i=0, \\
&pM_f,\ \ & \text{ if } f(1)\in \Ib,\ i=0.
\end{aligned}
\right.
\]

A weight $f\in \Ibw^d$ is called {\em anti-dominant} if
$\frac{m-1}{2}\geq f(1)\geq f(2)\geq \cdots \geq f(d);$ 
in this case we have $f(j) \in \Iwl \cup \Ib$, for all $j$. 
Denote $\Ianti=\{f\in \Ibw^d \mid f \text{ is anti-dominant} \}.$
Decompose $\V^{\otimes d}$ into a direct sum of cyclic submodules generated by $M_f$, for  anti-dominant weights $f$: $\V^{\otimes d}=\bigoplus_{f \in \Ianti} \M_f,$ where $\M_f =M_f\HB.$ Denote by $\mathcal O_f$ the orbit of $f$ under the action of $W_d$ on $\Ibw^d$. 
The $\HB$-module $\M_f$ admits a standard basis $\{M_g \mid g \in \mathcal O_f \}$.  

The stabilizer subgroup of $f \in \Ianti$ in $W_d$ is of the form
\begin{align*}
W_f =W_{{m_1}}\times \ldots \times W_{{m_k}}\times S_{m_{k+1}}\times \ldots\times S_{m_l},\end{align*} 
with all $m_i>0$ and $W_{{m_1}}\times \ldots \times W_{{m_k}}$ corresponding to the components of $f$ in $\Ib$. 
Note the stabilizer subgroup $W_f$ is not a parabolic subgroup of $W_d$ when $k\ge 2$. This phenomenon does not occur in the setting of \cite{BW18a, BWW20}. 
We call the summands $\M_f$ of $\V^{\otimes d}$ {\em quasi-permutation modules}. Clearly, for $f,f' \in \Ianti$, we have $\M_f \cong \M_{f'},$ if $W_f =W_{f'}.$  
If $W_f$ is not parabolic, $\M_f$ is in general not an induced module as those considered in parabolic KL setting \cite{De87}.

Let $f \in \Ianti$. Denote by ${}^fW$ the set of minimal length right coset representatives in $W_f \backslash W_d$. We define a $\Q$-linear map $\iba$ on the module $\M_f$ (which has a basis $M_{f\cdot \sigma}$, for $\sigma \in {}^fW$) by
$\iba(v)=v^{-1}, \iba(M_{f \cdot \sigma})=M_f\bar H_{\sigma}, \forall \sigma \in {}^fW.$
It can be shown \cite{SW21} (more difficult than the parabolic case in \cite{De87}) that the map $\iba$ on $\M_f$ is compatible with the bar operator on the Hecke algebra, i.e., 
$\iba(xh)=\iba(x) \overline{h},$ for all $x \in \M_f, \  h \in \HB$.
In particular, $\iba^2 =\text{Id}$. We shall call $\iba$ the bar involution on $\M_f$.

\begin{theorem}
\cite{SW21}
  \label{thm:CBMf}
Suppose $p \in v^\Z$, and let $f\in \Ianti$. Then for each $\sigma \in {}^fW$, there exists a unique element $C_{ \sigma} \in \M_f$ such that
\[
\iba (C_{\sigma}) =C_{\sigma},\;  \text{ and  }  \;
C_{\sigma} \in M_{f\cdot \sigma}+\sum_{w\in {}^fW, w< \sigma }\limits 
v^{-1}\Z[v^{-1}] M_{f\cdot w}.
\]
\end{theorem}

Similarly, there exist elements $C^*_{\sigma} \in \M_f$, for $\sigma \in {}^fW$, characterized by 
$\iba (C^*_{\sigma}) =C^*_{\sigma}$ and 
$C^*_{\sigma} \in M_{f\cdot \sigma}+\sum_{w\in {}^fW, w<\sigma} 
v \Z[v] M_{f\cdot w}.$
The basis $\{ C_{\sigma} \mid \sigma \in {}^fW \}$ is called a {\em quasi-parabolic KL basis} for $\M_f$; the basis $\{ C^*_{\sigma} \mid \sigma \in {}^fW \}$ is called a {\em dual quasi-parabolic KL basis} for $\M_f$.  Depending on the choice of $f$, the canonical basis can be type B or type A parabolic KL basis \cite{KL79, De87}, or neither. 

\subsection{A type B q-Schur duality}
  \label{sec:Schur}

Set $n =\frac{m}2 \in \frac12 \N.$ We consider the quantum group $\bU = \bU_v(\mathfrak {sl}_{2r+m})$ of type $A_{2r+m-1}$, where $\I :=  [1-n -r .. n+r-1 ].$ We view $\V$ as a natural representation of $\bU$, and so $\V^{\otimes d}$ is a $\bU$-module via the comultiplication $\Delta$. 
We consider the Satake diagram of type AIII with $m-1 =2n-1$ black nodes and $r$ pairs of white nodes, and a diagram involution $\tau$:
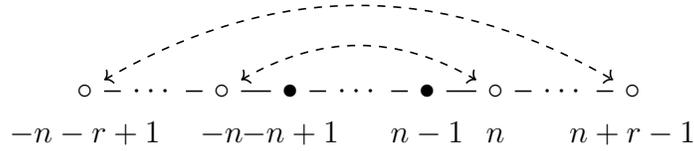
\begin{center}
\begin{tikzpicture}[scale=1, semithick]
\node (-4) [label=below:{$-n-r+1$}] at (0,0){$\circ$};
\node (-3)  at (1.3*0.7,0) {$\cdots$} ;
\node (-2) [label=below:{$-n$}] at (2.6*0.7,0){$\circ$};
\node (-1) [label=below:{$-n+1$}] at (3.9*0.7,0){$\bullet$};
\node (0)  at (5.2*0.7,0){$\cdots$};
\node (1) [label=below:{$n-1$}] at (6.5*0.7,0){$\bullet$};
\node (2) [label=below:{$n$}] at (7.8*0.7,0){$\circ$};
\node (3)  at (9.1*0.7,0){$\cdots$};
\node (4) [label=below:{$n+r-1$}] at (10.4*0.7,0){$\circ$};
\path (-4) edge (-3)
          (-3) edge (-2)
          (-2) edge (-1)
          (-1) edge (0)
          (0) edge (1)
          (1) edge (2)
          (2) edge (3)
          (3) edge (4);
\path (-4) edge[dashed,bend left,<->] (4)
          (-2) edge[dashed,bend left,<->] (2);
\end{tikzpicture}
\end{center}
(In case $n=0$, the black nodes are dropped; the nodes $n$ and $-n$ are identified and fixed by $\tau$.) The involution $\tau$ on $\I$ sends $i \mapsto \tau(i)= -i$, for all $i$, and it induces an involution of $\bU$, denoted again by $\tau$, by permuting the indices of its generators $E_i, F_i, K_i^{\pm 1}$. 

Let $\I_\bullet =[1-n.. n-1]$ be the set of all black nodes in $\I$ and $\I_\circ = \I \backslash \I_\bullet$. Associated to the Satake diagram $(\I =\I_\bullet \cup \I_\circ, \tau)$, we have a quantum symmetric pair $(\bU, \Ui)$ of type AIII. Recall that $p \in v^\Z$. 
We shall fix the parameters to be 
 \begin{align}
   \label{eq:Bi0}
\left\{ 
\begin{aligned}
\va_i &=1\ (\text{for } i\neq \pm n), \\
 \va_{n} &=(-1)^{m-1}v^{m}p^{-1},  \quad
 \va_{-n} =p,  \qquad\qquad   \text{ if } m=2n \in \Z_{\ge 1};
 \end{aligned}
\right.
\end{align}
 \begin{align}
   \label{eq:Bi1}
\left\{ 
\begin{aligned}
\va_{0} & =v^{-1}, \quad
\va_i =1\ (\text{for } i\neq 0),   \\
 \kappa_0 &= (p-p^{-1}) / (v-v^{-1}),
 \qquad\qquad\qquad\qquad \text{ if } m=0.
 \end{aligned}
\right.
\end{align}

It can be shown \cite{SW21}
that the actions of $\Ui$ and $\HB$ on $\V^{\otimes d}$ commute with each other: 
\[
\Ui \stackrel{\Psi}{\curvearrowright} \V^{\otimes d} \stackrel{\Phi}{\curvearrowleft} \Hy_{B_d}. 
\]
Moreover, $\Psi(\Ui)$ and $\Phi(\HB)$ form double centralizers in $\text{End} (\V^{\otimes d})$. This duality has been called {\em $\imath$Schur duality} (which goes back to \cite{BW18a} in a quasi-split case). 

There exists a unique anti-linear bar involution
$ \iba\colon \V^{\otimes d} \rightarrow \V^{\otimes d}$ such that $\iba(M_f) =M_f$, for $f\in \Ianti$, and it is compatible with the bar involutions on $\Hy_{B_d}$ and $\Ui$; that is, for $u\in \Ui$, $v\in \V^{\otimes d}$, and $h \in \HB$,  
$\iba(uvh)=\iba(u)\iba(v)\bar h.$
Recall that $\V^{\otimes d}$ is a direct sum of the quasi-permutation modules $\M_f$ of $\HB$. The union of the (dual) quasi-parabolic KL bases on the summands $\M_f$ provides a (dual) quasi-parabolic KL basis on $\V^{\otimes d}$.

\begin{theorem} 
\cite{SW21}
  \label{thm:iCBsame}
The (dual) $\imath$canonical basis on $\V^{\otimes d}$ (viewed as a $\Ui$-module) coincides with the (dual) quasi-parabolic KL basis on $\V^{\otimes d} =\bigoplus_{f} \M_f$ (viewed as an $\HB$-module). 
\end{theorem}

The quasi-parabolic KL polynomials are by definition the transition matrix entries from the quasi-parabolic KL to the standard basis. An inversion formula for parabolic KL polynomials (theorems of Kazhdan-Lusztig and Douglass) can be generalized to the quasi-parabolic cases. 

\begin{remark}
In case when $r=0$, $\imath$Schur duality reduces to Jimbo-Schur duality \cite{Jim86} between $\bU$ and Hecke algebra of type A, and Theorem~\ref{thm:iCBsame} goes back to a result of Frenkel-Khovanov-Kirillov. In case when $m=0, 1$, it reduces to the quasi-split $\imath$Schur duality \cite{BW18a} \cite{Bao17} \cite{BWW18}, which has applications to Kazhdan-Lusztig theory of classical type. 
\end{remark}

\section{Application to super Kazhdan-Lusztig theory}

\subsection{The BGG category}

Consider the BGG category $\mathcal{O}$ of $\mathfrak{g}$-modules, where $\mathfrak{g}=\mf{n}^-\oplus\mf{h}\oplus\mf{n}$ is a simple or reductive Lie (super) algebra over $\mathbb{C}$. 
There is a duality functor $\vee:\mathcal{O}\rightarrow\mathcal{O}$ sending $M=\oplus_{\mu\in\mathfrak{h}^*}M_\mu$ to $M^\vee:=\oplus_{\mu\in\mathfrak{h}^*}M_\mu^*$. 
Let $M(\lambda)$ be the Verma module with highest weight $\lambda$ and $L(\lambda)$ be its unique irreducible quotient. 
It is known that a simple module $L(\lambda)$ and a tilting module $T(\lambda)$ (of highest weight $\lambda$) are self-dual with respect to $\vee$.

For semisimple Lie algebras, the linkage is controlled by the dot action of the Weyl group, and the BGG category $\mathcal O$ admits a block decomposition according to the central characters. For a (or any) regular block $\mc O^0$, its Grothendieck group is identified with the Weyl group algebra. If we further identify the Verma module basis in $[\mc O^0]$ with the standard basis in the Hecke algebra (specialized at $v=1$), then the Kazhdan-Lusztig conjecture (a theorem of Beilinson-Bernstein and Brylinski-Kashiwara) states that the simple module basis corresponds to the dual canonical basis (and the tilting module basis corresponds to the canonical basis). 

For general linear or ortho-symplectic Lie superalgebras, the linkage in the BGG category is no longer controlled by the Weyl group, and so the formulation of Kazhdan-Lusztig theory via Hecke algebras breaks down.

\subsection{Super type BCD character formulas}

Let us treat the super type B case, $\mathfrak{g}=\mathfrak{osp}(2m+1|2n)$, in detail. With respect to a standard Dynkin diagram
\begin{center}
\vspace{-5mm}
\hskip -3cm \setlength{\unitlength}{0.16in}
\begin{picture}(24,4)
\put(5.6,2){\makebox(0,0)[c]{$\bigcirc$}}
\put(8,2){\makebox(0,0)[c]{$\bigcirc$}}
\put(10.4,2){\makebox(0,0)[c]{$\bigcirc$}}
\put(14.85,2){\makebox(0,0)[c]{$\bigotimes$}}
\put(17.25,2){\makebox(0,0)[c]{$\bigcirc$}}
\put(19.4,2){\makebox(0,0)[c]{$\bigcirc$}}
\put(23.5,2){\makebox(0,0)[c]{$\bigcirc$}}
\put(8.35,2){\line(1,0){1.5}} \put(10.82,2){\line(1,0){0.8}}
\put(13.2,2){\line(1,0){1.2}} \put(15.28,2){\line(1,0){1.45}}
\put(17.7,2){\line(1,0){1.25}} \put(19.81,2){\line(1,0){0.9}}
\put(22,2){\line(1,0){1}}
\put(6.8,2){\makebox(0,0)[c]{$\Longleftarrow$}}
\put(12.5,1.95){\makebox(0,0)[c]{$\cdots$}}
\put(21.5,1.95){\makebox(0,0)[c]{$\cdots$}}
\put(5.4,1){\makebox(0,0)[c]{\tiny $\epsilon_{1}$}}
\put(7.8,1){\makebox(0,0)[c]{\tiny $-\epsilon_{1} +\epsilon_{2}$}}
\put(14.4,1){\makebox(0,0)[c]{\tiny $-\epsilon_{m} +\delta_{1}$}}
\put(17.3,1){\makebox(0,0)[c]{\tiny $-\delta_{1} +\delta_{2}$}}
\put(23.5,1){\makebox(0,0)[c]{\tiny $-\delta_{n-1} +\delta_{n}$}}
\end{picture}
\vspace{-3mm}
\end{center}
we have the Weyl vector $\scriptstyle \rho = \frac12 \epsilon_1 +\frac32 \epsilon_2 +\ldots  +(m-\frac12) \epsilon_{m}  -(m-\frac12) \delta_1 - (m-\frac32) \delta_2 -\ldots - (m-n +\frac12) \delta_n.$ 
There exists a $\rho$-shift bijection for the set of integer weights
\begin{align*}
X^{m|n} :=\oplus_{i=1}^m\Z\epsilon_i\oplus\oplus_{j=1}^n \Z\delta_j
  \xrightarrow[\mbox{\tiny $\rho$-shift}]{\cong}({\textstyle \frac{1}{2} } +\Z)^{m+n}, \quad \lambda\mapsto f_\lambda,
  \end{align*}
 where $f_\lambda$ is defined via $\lambda+\rho =\sum_{i=1}^m f_\lambda(i) \epsilon_i + \sum_{j=1}^n f_\lambda(m+j) \delta_j$.
Similarly, there exists a bijection for the set of half integer weights
$
  X^{m|n}_{{\scriptstyle \frac{1}{2}}}
  :=\oplus_{i=1}^m({\textstyle \frac{1}{2} } +\Z)\epsilon_i\bigoplus
  \oplus_{j=1}^n ({\textstyle \frac{1}{2} } +\Z)\delta_j
  \xrightarrow[\mbox{\tiny $\rho$-shift}]{\cong}\Z^{m+n},$ $\lambda\mapsto f_\lambda.
$ 
Denote by $\mathcal{O}^{m|n}_{\mf b}$ (respectively, $\mathcal{O}^{m|n}_{\mf b, \frac12}$) the BGG category which contains the Verma modules $M(\lambda)$, tilting modules $T(\lambda)$ and simple modules $L(\lambda)$, parametrized by the weights $\la \in X^{m|n}$ (respectively, $\la \in X^{m|n}_{{\scriptstyle \frac{1}{2}}}$). 

Recall from \S\ref{sec:Schur} the quasi-split quantum symmetric pair $(\bU_v(\mathfrak{sl}_N), \Ui)$ of type AIII, where we fix $p=v$ in \eqref{eq:Bi0}-\eqref{eq:Bi1}. Recall the natural representation $\V$ with basis $\{u_i\mid i \in [\textstyle{\frac{1-N}{2}..\frac{N-1}{2}}]\}$, for $N$ even and odd, allowing $N =\infty$ with parity! Then we can identify the indexing set $[\textstyle{\frac{1-N}{2}..\frac{N-1}{2}}]$ with ${\textstyle \frac{1}{2} } +\Z$  for $N=\infty$ (even), and with $\Z$  for $N=\infty$ (odd).

  By Theorem~\ref{thm:iCBmodule}, the $\bU_v(\mathfrak{sl}_{\infty})$-module $\V^{\otimes m}\otimes \V^{*\otimes n}$ (regarded as $\Ui$-module with $p=v$) admits an $\imath$canonical basis, denoted by $\{C^\imath_f\}$, and a dual $\imath$canonical basis, denoted by $\{L^\imath_f\}$, where $f \in ({\textstyle \frac{1}{2} } +\Z)^{m+n}$ or $\Z^{m+n}$, respectively.

  Define the following $\Z$-module isomorphisms
  \begin{align}
  \label{eq:P2}
  \begin{split}
  \Psi_{\mf b}: [\mathcal{O}^{m|n}_{\mf b}]\longrightarrow \V_\Z^{\otimes m} \otimes \V_{\Z}^{*\otimes n}, & \quad
  [M(\lambda)]\mapsto M_{f_\lambda}\;\; (\la \in X^{m|n}),
  \\
  \Psi_{\mf b, \frac12}: [\mathcal{O}^{m|n}_{\mf b,\frac12}]\longrightarrow \V_\Z^{\otimes m} \otimes \V_{\Z}^{*\otimes n}, & \quad
  [M(\lambda)]\mapsto M_{f_\lambda} \;\; (\la \in X^{m|n}_{{\scriptstyle \frac12}}).
  \end{split}
  \end{align}
  A basic fact here \cite{BW18a} is that the generators $B_i$ in $\Ui$ act on $[\mathcal{O}^{m|n}_{\mf b}]$ and $[\mathcal{O}^{m|n}_{\mf b,\frac12}]$ by translation functors, and the above $\Z$-module isomorphisms become $\Ui_\Z$-module isomorphisms at $v=1$. (A similar observation on translation functors and $B_i$ is valid for $p=1$ \cite{Bao17}, and it was made independently in \cite{ES18} in  the non-super setting.)

\begin{theorem} \cite{BW18a}
  \label{thm:KLb}
  The $\Z$-module isomorphism $  \Psi_{\mf b}$ (respectively, $\Psi_{\mf b, \frac12}$) in \eqref{eq:P2} sends
  \begin{equation*}
  [L(\lambda)]\mapsto L^\imath_{f_{\lambda}}, \quad [T(\lambda)]\mapsto C^\imath_{f_{\lambda}},
  \quad
 \text{ for $\la \in X^{m|n}$ (and respectively, $\la \in X^{m|n}_{{\scriptstyle \frac12}}$).}
 \end{equation*}
\end{theorem}

\begin{remark} [Super type A Kazhdan-Lusztig theory]
Consider the BGG category $\mathcal{O}^{m|n}$ of modules over the general linear Lie superalgebra $\mathfrak{g}=\mathfrak{gl}(m|n)$ of integer weights. We have an almost identical $\Z$-module isomorphism $\Psi: [\mathcal{O}^{m|n}] \rightarrow \V_\Z^{\otimes m} \otimes \V_{\Z}^{*\otimes n}$ as in \eqref{eq:P2}, which match the Verma basis with the standard basis. Then Brundan-Kazhdan-Lusztig conjecture \cite{Bru03} (proved by Cheng, Lam, and the author in \cite{CLW15}) states that the simple module basis (respectively, tilting module basis) is mapped by $\Psi$ to Lusztig dual canonical basis (respectively, canonical basis). There has been a second proof in \cite{BLW17} using ideas of categorification.
\end{remark}

\begin{example}
  Take $n=0$ and $m=1$, so that $\mathfrak{g}=\mathfrak{so}_3\cong\mathfrak{sl}_2$. If the standard basis (= canonical basis) $\{u_i \mid i\in\Z\}$ for $\V$ is indexed by $\Z$, then $\V$ admits an $\imath$canonical basis $\{u_0, u_{-i}, u_{i}+ v^{-1} u_{-i} \mid i\in\Z_{>0}\}$ and a dual $\imath$canonical basis $\{u_0, u_{-i},u_{i}- v u_{-i} \mid i\in\Z_{>0}\}$. 
\end{example}

Theorem~\ref{thm:KLb} can be adapted to the type D Lie superalgebra $\mathfrak{g}=\mathfrak{osp}(2m|2n)$, by setting the parameter $p=1$ (instead of $p=v$) in \eqref{eq:Bi0}-\eqref{eq:Bi1} (see Bao \cite{Bao17}). Thanks to Theorem~\ref{thm:iCBsame}, Theorem~\ref{thm:KLb} for $m=0$ amounts to a reformulation for the type B Kazhdan-Luszitg conjecture \cite{KL79}. For further extension to Kazhdan-Lusztig theory for super parabolic BGG categories, see \cite{BWW20}. 

\section{Hall algebras}

The Drinfeld double quantum group $\tU=\tU(\g)$ is the $\Qq$-algebra generated by $E_i,F_i,K_i,{K_i'}$ $(i\in \I)$ subject to relations in \cite[(6.1)-(6.5)]{LW19a} similar to those in $\bU$. Denote the Cartan matrix for $\g$ by $(c_{ij})_{i,j\in \I}$. 
Following \cite{LW19a}, we define the \emph{universal $\imath$quantum group} $\tUi$ associated to a Satake diagram $(\I=\Iblack \cup \Iwhite,\tau)$ as the $\Qq$-subalgebra of $\tU$ generated by  $\tU_{\I_\bullet}$ and $\{B_i,\tk_i\mid i\in \Iwhite \}$, with identifications
\begin{align}
B_i&\mapsto  F_i + \tT_{w_\bullet}(E_{\tau i}) K'_i,\qquad \tk_i \mapsto K_i K'_{\tau i}, \qquad i\in \Iwhite.
\end{align}
Denote the embedding by $\imath: \tUi \rightarrow \tU$, $x \mapsto x^\imath$. One checks that $\tUi$ is a right coideal subalgebra of $\tU$. 
The $\imath$quantum group $\Ui$ (with $Y =\N \I$) can be obtained from $\tUi$ via a central reduction.

In the remainder of this section, we shall only consider $\tUi$ of quasi-split types, i.e., $\Iblack =\emptyset$. Let $Q=(Q_0,Q_1)$ be a virtually acyclic quiver; see \cite[Def.~4.4]{LW20a}. This is a mild generalization of acyclic quivers, allowing a generalized Kronecker subquiver. Throughout the paper, we shall identify $Q_0=\I$. An {\em $\imath$quiver} $(Q, \tau)$ consists of a (virtually acyclic) quiver $Q$ and an involution $\tau$ of $Q$; we allow the {\em trivial} involution $\Id$. We work over a finite field $\mathbb F_q$. 
An involution $\tau$ of $Q$ induces an involution of the path algebra $\mathbb F_q Q$, also denoted by $\tau$.

Let $\ov{Q}$ be a new quiver obtained from $Q$ by adding a loop $\varepsilon_i$ at the vertex $i\in Q_0$ if $\tau i=i$, and adding an arrow $\varepsilon_i: i\rightarrow \tau i$ for each $i\in Q_0$ if $\tau i\neq i$; the $\varepsilon_i$ are in purple color below. The \emph{$\imath$quiver algebra} $\Lambda^\imath$ associated to $(Q, \tau)$ can be defined in terms of the quiver $\ov Q$ with relations, cf. \cite[Proposition 2.6]{LW19a}; that is, $\Lambda^\imath \cong \mathbb F_q \ov{Q} / \ov{I}$, where $\ov{I}$ is generated by
$\varepsilon_{i}\varepsilon_{\tau i}$ for each $i\in Q_0$ and 
$\varepsilon_i \alpha -\tau(\alpha)\varepsilon_j$ for each arrow $\alpha:j\rightarrow i$ in $Q_1$.

Rank 1 or 2 subquivers of the quiver $\ov{Q}$ associated to a general virtually acyclic quiver $Q$ look like as follows (where $\ov{Q}$ is obtained from $Q$ by adding arrows $\purple{\varepsilon}$'s): 
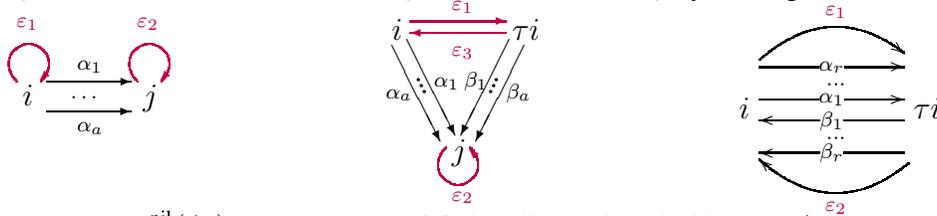
\begin{center}\setlength{\unitlength}{0.7mm}
\vspace{-12mm}
\begin{equation*}
\begin{picture}(100,50)(0,10)
{\setlength{\unitlength}{0.8mm}
\put(-30,30){\begin{picture}(50,3)(0,0) \put(0,-3){$i$}
\put(4,1){\vector(1,0){14}}
\put(4,-4){\vector(1,0){14}}
\put(8,-2.5){\tiny$\cdots$}
\put(9,1){$^{\alpha_1}$}
\put(9,-6.5){$_{\alpha_a}$}
\put(20,-3){$j$}
\color{purple}
\qbezier(-1,1)(-3,3)(-2,5.5)
\qbezier(-2,5.5)(1,9)(4,5.5)
\qbezier(4,5.5)(5,3)(3,1)
\put(3.1,1.4){\vector(-1,-1){0.3}}
\qbezier(19,1)(17,3)(18,5.5)
\qbezier(18,5.5)(21,9)(24,5.5)
\qbezier(24,5.5)(25,3)(23,1)
\put(23.1,1.4){\vector(-1,-1){0.3}}
\put(-1,10){\tiny$\varepsilon_1$}
\put(19,10){\tiny$\varepsilon_2$}
\end{picture}}
}
{\setlength{\unitlength}{0.8mm}
\put(30,30){ \begin{picture}(50,10)(0,-10)
\put(0,-2){$i$}
\put(20,-2){$\tau i$}
\put(-1,-11){$_{\alpha_a}$}
\put(7,-9){$_{\alpha_1}$}
\put(12,-9){$_{\beta_1}$}
\put(19,-11){$_{\beta_a}$}
\put(3.7,-10){$\cdot$}
\put(4.2,-11){$\cdot$}
\put(4.7,-12){$\cdot$}

\put(0,-3){\vector(1,-2){8}}
\put(2.5,-2){\vector(1,-2){8}}
\put(19.5,-2){\vector(-1,-2){8}}
\put(16.7,-10){$\cdot$}
\put(16.2,-11){$\cdot$}
\put(15.7,-12){$\cdot$}
\put(22,-3){\vector(-1,-2){8}}
\put(10,-22){$j$}
\color{purple}
\put(3,1){\vector(1,0){16}}
\put(19,-1){\vector(-1,0){16}}
\put(10,1){$^{\varepsilon_1}$}
\put(10,-4){$_{\varepsilon_3}$}
\put(10,-28){$_{\varepsilon_2}$}
\begin{picture}(50,23)(-10,19)
\color{purple}
\qbezier(-1,-1)(-3,-3)(-2,-5.5)
\qbezier(-2,-5.5)(1,-9)(4,-5.5)
\qbezier(4,-5.5)(5,-3)(3,-1)
\put(3.1,-1.4){\vector(-1,1){0.3}}
\end{picture}
\end{picture} }
}
\put(110,30){\begin{picture}(50,43)(10,0)
\xymatrix{ {i }\ar@/^1.5pc/@<2.5ex>@[purple][rr]^{\textcolor{purple}{\varepsilon_1}} \ar@<3ex>[rr]|-{\alpha_r}   \ar@<0.75ex>[rr]|-{\alpha_1}^{\cdots}   
&& {
\tau i}\ar@<0.75ex>[ll]|-{\beta_1} \ar@/^1.5pc/@<2.5ex>@[purple][ll]^{\textcolor{purple}{\varepsilon_2}} \ar@<3ex>[ll]|-{\beta_r}_{\cdots}  }
\end{picture} }
\end{picture}
\vspace{-2mm}
\end{equation*}
\end{center}

Denote by $\text{mod}^{\text{nil}}(\Lambda^\imath)$ the category of finite-dimensional nilpotent $\Lambda^\imath$-modules. Denote by $S_i$ the 1-dimensional $\Lambda^\imath$-module supported at $i \in \I$, and $\K_i$ the 2-dimensional module ``supported at $\varepsilon_i$".
The algebra $\Lambda^\imath$ is a 1-Gorenstein algebra and hence admits favorable homological properties. In particular, the subcategory $\cp^{\le 1}(\Lambda^\imath)$ of $\text{mod}^{\text{nil}}(\Lambda^\imath)$ consisting of modules of projective dimension at most $1$ admits clean characterization.

Let $\mathcal H (\Lambda^\imath)$ be the Ringel-Hall algebra of $\text{mod}^{\text{nil}}(\Lambda^\imath)$ over $\Q(\sqrt{q})$, that is,  the $\Q(\sqrt{q})$-vector space whose basis is formed by the isoclasses $[M]$ of objects $M\in \text{mod}^{\text{nil}}(\Lambda^\imath)$, with multiplication defined by  
$[M]\diamond [N]=\sum_{[L]\in \text{Iso}(\text{mod}(\Lambda^\imath))}\frac{|\text{Ext}^1(M,N)_L|}{|\Hom(M,N)|}[L].$
Then, the semi-derived Hall algebra $\mathcal{SDH}(\Lambda^\imath)$ of $\Lambda^\imath$ is defined in terms of localization of a quotient algebra of the Ringel-Hall algebra $\mathcal H (\Lambda^\imath)$ with respect to $\cp^{\le 1}(\Lambda^\imath)$, and the \emph{$\imath$Hall algebra} $\tMHk$ is defined to be $\mathcal{SDH}(\Lambda^\imath)$ with a new multiplication via  twisting by an Euler form; see \cite[Appendix~A]{LW19a} \cite{LW20a} for precise definitions. 

Let $\I_\tau$ be a set of representatives of the $\tau$-orbits on $\Iwhite$. 

\begin{theorem}
\cite{LW19a, LW20a}
   \label{thm:main}
Let $(Q, \tau)$ be a virtually acyclic $\imath$quiver. Then there exists a $\Q(\sqrt{q})$-algebra monomorphism $\widetilde{\psi}: \tUi_{|v= \sqrt{q}}  \rightarrow \tMHk,$ which sends
\begin{align*}
B_j \mapsto \frac{-1}{q-1}[S_{j}],\text{ if } j\in\I_\tau,
&\qquad\qquad
\tk_i \mapsto - q^{-1}[\bK_i], \text{ if }\tau i=i \in \I,
\\
B_{j} \mapsto \frac{{\sqrt{q}}}{q-1}[S_{j}],\text{ if }j\notin \I_\tau,
&\qquad\qquad
\tk_i \mapsto \sqrt{q}^{\frac{-c_{i,\tau i}}{2}}[\bK_i],\quad \text{ if }\tau i\neq i \in \I.
 \end{align*}
\end{theorem}
This theorem for diagonal $\imath$quivers  $(Q \sqcup Q, \text{swap})$ specializes to 
 Bridgeland's Hall algebra realization of Drinfeld double quantum groups in \cite{Br13}. 
The above monomorphism becomes an isomorphism for Dynkin $\imath$quivers. 
Reflection functors on $\imath$Hall algebras provide a conceptual approach to braid group actions on $\tUi$ (of quasi-split type); see \cite{LW21a, LW21b}.

\section{Relative braid group actions}
\label{sec:braid}

Let $(\tU,\tUi)$ be the universal quantum symmetric pair associated to a Satake diagram $(\I=\Iblack \cup \Iwhite,\tau)$. It is shown in \cite{WZ21} that there exists a quasi $K$-matrix $\tfX$ associated to $(\tU,\tUi)$ satisfying an intertwining relation like \eqref{eq:Upsilon}. 
%
For $i\in \Iwhite$, denote by $\tfX_i$ the quasi $K$-matrix associated to the rank one Satake subdiagram $\I_{\bullet,i}:=\Iblack \cup \{i,\tau i\}$ in the setting of $(\tU,\tUi)$. 
 
 Let $W_{\bullet,i}$ be the parabolic subgroup of the Weyl group $W$ of $\g$ generated by $s_i$, for $i\in \I_{\bullet,i}$, with longest element $w_{\bullet,i}$. Define $\bs_i \in W_{\bullet,i}$ by $\bs_i  =\bwi w_\bullet$; it is clear that $\bs_i=\bs_{\tau i}$. Recall $\I_\tau$ denotes a set of representatives of the $\tau$-orbits on $\Iwhite$. 
The relative Weyl group associated to the symmetric pair is identified with the subgroup $\bbW=\langle \bs_i\mid i\in \I_\tau \rangle$ of $W$. There is a notion of the relative braid group associated to $\bbW$. The existence of such a relative braid group action on an $\imath$quantum group $\Ui$ was conjectured in \cite{KP11}. The conjecture was verified therein via compute computation for (mostly quasi-split) finite type; for an alternative approach via $\imath$Hall algebras see \cite{LW21a, LW21b}. 

There is a braid group action associated to $W$ on the Drinfeld double $\tU$ (see \cite{LW21b}), a variant of the braid group action on $\bU$ in \cite{Lus93}. We shall need a suitably rescaled variant, denote by $\tcT_i^{-1}$, for $i\in \I$, which again satisfies the braid group relations. In particular, an automorphism $\tcT_{\bs_i}^{-1}$ of $\tU$, for $i\in \I_\tau$, is defined. We announce a new conceptual approach developed with W.~Zhang to relative braid group actions. 

\begin{theorem} \cite{WZ21} 
  \label{thm:newb0}
Let $(\tU, \tUi)$ be a universal quantum symmetric pair of arbitrary finite type.  
Then there exists an automorphism $\tTT_i^{-1}$ of $\tUi$, for $i\in \I_\tau$, which satisfies the intertwining relation
\begin{align*}
\tTT_i^{-1}(x)^\imath \cdot \tfX_i =\tfX_i \cdot \tcT_{\bs_i}^{-1}(x^\imath),
\quad \text{ for all } x\in \tUi. 
\end{align*}
Moreover, the automorphisms $\tTT_i^{-1}$, for $i\in \I_\tau$, satisfy the relative braid group relations. 
\end{theorem}
The approach in \cite{WZ21} has additional consequences. Explicit compact formulas for the action of $\tTT_i^{-1}$ on the generators of $\tUi$ are obtained. 
The relative braid group action on $\tUi$ gives rise to compatible relative braid group actions on $\Ui$  and $\bU$-modules (viewed as $\Ui$-modules). Along the way, we prove the conjecture of Dobson-Kolb \cite{DK19} on factorization of quasi K-matrices of arbitrary finite type. 

\section{A current presentation of affine type}
  \label{sec:affine}

In this section, we consider universal $\imath$quantum groups $\tUi$ of split affine type, that is, $\Iblack =\emptyset, \tau =\Id$, and the Cartan matrix $(c_{ij})_{i,j\in \I}$ is of untwisted affine type. By definition, $\tUi$ is a subalgebra of $\tU$; alternatively, $\tUi$ is the $\Qq$-algebra generated by $B_i, \K_i^{\pm 1}$ $(i\in \I)$, subject to the following relations: $\K_i$  are central, and 
\begin{align}
B_iB_j -B_j B_i&=0, \quad \qquad\qquad\qquad\qquad\qquad \text{ if } c_{ij}=0,
 \label{eq:S1} \\
B_i^2 B_j -[2]_{i} B_i B_j B_i +B_j B_i^2 &= - v_i^{-1}  B_j \K_i,  \quad\qquad\qquad\qquad \text{ if }c_{ij}=-1.
 \label{eq:S2}
\end{align}
We omit here the more complicated Serre type relations between $B_i, B_j$ for $c_{ij}=-2, -3$;  they can be read off from setting $k=l=0$ in \eqref{iDRG5'}--\eqref{iDRG6'} below. The $\K_i$ (which is natural from $\imath$Hall algebra viewpoint) is related to $\tk_i$ used earlier by $\K_i = -v_i^2 \tk_i$. 

Associated to the affine Lie algebra $\g$, we denote by $\g_0, \I_0, X_0, W_0$ the underlying semisimple Lie algebra, its simple roots, weight lattice and finite Weyl group. Recall the (extended) affine Weyl group $W^e =W_0 \ltimes X_0$. There are automorphisms $\tTT_i$ of $\tUi$, for $i\in \I$, which arise naturally from $\imath$Hall algebras \cite{LW21a, LW21b}; also see \S\ref{sec:braid}. They give rise to the action of an affine braid group associated to $W^e$ for the affine Lie algebra $\g$. In particular, we have automorphisms $\tTT_w$ of $\tUi$, for $w\in W^e$.

Define a sign function
$o(\cdot): \I_0 \longrightarrow \{\pm 1\}$
such that $o(i) o(j)=-1$ whenever $c_{ij} <0$. 
Define $v$-root vectors $B_{i,k},\acute{\Theta}_{i,m},\Theta_{i,m}$ in $\tUi$ for $i\in \I_0 $, $k\in \Z$ and $m\ge 1$ by \cite{LW20b, Zh21} 
\begin{align*}
B_{i,k} &= o(i)^k \tTT_{\omega_i}^{-k} (B_i),
  \\
\acute{\Theta}_{i,m} &=  o(i)^m \Big(-B_{i,m-1} \tTT_{\omega_i'} (B_i) +v_i^{2} \tTT_{\omega_i'} (B_i) B_{i,m-1}
+ (v_i^{2}-1)\sum_{p=0}^{m-2} B_{i,p} B_{i,m-p-2}  \K_{i}^{-1}\K_{\delta} \Big),
\notag \\
\Theta_{i,m} &=\acute{\Theta}_{i,m} - \sum\limits_{a=1}^{\lfloor\frac{m-1}{2}\rfloor}(v_i^2-1) v_i^{-2a} \acute{\Theta}_{i,m-2a}\K_{a\delta} -\delta_{m,ev} v_i^{1-m} \K_{\frac{m}{2}\delta}.
\end{align*}
A version of the $v$-root vectors for $\Ui$ in affine rank 1 case (known as q-Onsager algebra) was constructed earlier in \cite{BK20}.

Let $ \tUiD$ be the $\Q(v)$-algebra generated by $\K_{i}^{\pm1}$, $C^{\pm1}$, $H_{i,m}$ and $B_{i,l}$, where  $i\in \I_0$, $m \in \Z_{\geq 1}$, $l\in\Z$, subject to the following relations, for $m,n \in \Z_{\geq1}$ and $k,l\in \Z$:
\begin{align}
& \K_i, C \text{ are central, }\quad
[H_{i,m},H_{j,n}]=0, \quad \K_i\K_i^{-1}=1, \;\; C C^{-1}=1,\label{iDRG0}
\\
&[H_{i,m},B_{j,l}]=\frac{[mc_{ij}]_{i}}{m} B_{j,l+m}-\frac{[mc_{ij}]_{i}}{m} B_{j,l-m}C^m,\label{iDRG1'}
\\
&[B_{i,k}, B_{j,l+1}]_{v_i^{-c_{ij}}}  -v_i^{-c_{ij}} [B_{i,k+1}, B_{j,l}]_{v_i^{c_{ij}}}=0, \text{ if }i\neq j,\label{iDRG2'}
 \\ %
&[B_{i,k}, B_{i,l+1}]_{v_i^{-2}}  -v_i^{-2} [B_{i,k+1}, B_{i,l}]_{v_i^{2}}
=v_i^{-2}\Theta_{i,l-k+1} C^k \K_i-v_i^{-4}\Theta_{i,l-k-1} C^{k+1} \K_i\label{iDRG3'}
\\
&\qquad\qquad\qquad\qquad\qquad\qquad\quad\quad\quad
  +v_i^{-2}\Theta_{i,k-l+1} C^l \K_i-v_i^{-4}\Theta_{i,k-l-1} C^{l+1} \K_i, \notag
\\\label{iDRG34}
&[B_{i,k} ,B_{j,l}]=0,\qquad   \text{ if }c_{ij}=0,
\\
&   \sum_{s=0}^2(-1)^s
\qbinom{2}{s}_{i}
B_{i,k}^{2-s}B_{j,l} B_{i,k}^s =-v_i^{-1}   B_{j,l} \K_i C^k,\qquad\text{ if } c_{ij}=-1,\label{iDRG4'}
\\
\label{iDRG5'}
&\sum_{s=0}^3(-1)^s
\qbinom{3}{s}_{i}
B_{i,k}^{3-s}B_{j,l} B_{i,k}^s =-v_i^{-1} [2]_{i}^2 (B_{i,k} B_{j,l}-B_{j,l} B_{i,k})\K_i C^k,\quad\text{ if } c_{ij}=-2,
\\\label{iDRG6'}
&\sum_{s=0}^4(-1)^s
\qbinom{4}{s}_{i}
B_{i,k}^{4-s}B_{j,l} B_{i,k}^s = -v_i^{-1}(1+[3]_{i}^2)( B_{j,l} B_{i,k}^2+ B_{i,k}^2 B_{j,l})\K_i C^k\\\notag
&\qquad  +v_i^{-1}[4]_{i} (1+[2]_{i}^2) B_{i,k} B_{j,l} B_{i,k} \K_i C^k-v_i^{-2}[3]^2_{i} B_{j,l} \K_i^2 C^{2k} ,\quad\text{ if } c_{ij}=-3,
\end{align}
where $\Theta_{i,m}$ are related to $H_{i,m}$ by  
\begin{align*}
1+ \sum_{m\geq 1} (v_i-v_i^{-1})\Theta_{i,m} u^m  = \exp\Big( (v_i-v_i^{-1}) \sum_{m\geq 1} H_{i,m} u^m \Big).
\end{align*}
(The $\tUiD$ is denoted by $\tUiD_{\text{red}}$ in \cite{Zh21}.)

Below is an $\imath$-analog of the Drinfeld presentation of affine quantum groups \cite{Dr87} (proved in \cite{Be94, Da12}).

\begin{theorem}  \cite{LW20b, Zh21}
 \label{DprADE}
There is a $\Q(v)$-algebra isomorphism $\Phi: \tUiD \to \tUi$, which sends
\begin{align*}
B_{i,k} \mapsto B_{i,k},\quad H_{i,m} \mapsto H_{i,m},\quad \Theta_{i,m}\mapsto \Theta_{i,m},\quad \K_i \mapsto \K_i, \quad C \mapsto \K_\delta,
\end{align*}
for $i\in \I_0,k\in \Z, m\geq 1$.
\end{theorem}

\begin{remark}
More involved Serre relations among $B_{i,k}, B_{j,l}, B_{i,k'}$ generalizing the relations \eqref{iDRG4'}--\eqref{iDRG5'} are available; see \cite{LW20b, Zh21}. They can be shown to be equivalent to \eqref{iDRG4'}--\eqref{iDRG5'}, when combined with other relations \eqref{iDRG0}--\eqref{iDRG34} above. 
\end{remark}

It is straightforward to pass the $v$-root vectors and Drineld type presentation of $\tUi$ to $\Ui$ with arbitrary parameters by central reduction. This current (or Drinfeld type) presentation will be extended beyond split types in a future work.

\section{Open problems}
 \label{sec:open}

\vspace{1em}

\begin{quote}

{\em \purple{``There's no use trying, one can't believe i... things."}}

\noindent{\em \purple{``Why, sometimes I've believed as many as six i... things before breakfast."}}
\begin{flushright}{ \purple{--- Alice in Wonderland} } 
\end{flushright}
\end{quote}
\vspace{1em}

The open problems in the following six (interconnected) directions on $\imath$quantum groups look most appealing to us. 
\begin{enumerate}
\item
\emph{Positivity of $\imath$canonical basis}

Positivity of $\imath$canonical basis holds in (affine) type AIII \cite{LiW18} \cite{FL+20}. 
We conjecture that the $\imath$canonical bases arising from the $\imath$quantum groups of (quasi-) split ADE type (with parameters suitably specified; see Example~\ref{ex:rank1}) exhibit various positivity properties. 
Recently Lusztig extends his earlier construction of total positivity to symmetric spaces in \cite{Lus21}. It will be interesting to strengthen this construction by connecting to $\imath$canonical basis (with positivity). 

\item
 \emph{$\imath$Quiver varieties and geometric realizations of $\imath$quantum groups}
 
Geometric realizations of quantum groups are obtained in \cite{BLM90, V98, Nak00, MO18}. The works \cite{BKLW, FL+20} can be regarded as the $\imath$-generalizations of \cite{BLM90}. Li \cite{Li19} provides an $\imath$-analog of some Nakajima quiver varieties. We observe however that the diagram involutions used loc. cit. are in line with Vogan diagrams instead of Satake diagrams. A fresh start is needed to construct general $\imath$quiver varieties (allowing Satake diagrams with black nodes and non-Dynkin types).  The geometric realization of $\imath$quantum groups \`a la \cite{V98, Nak00, MO18} remains to be carried out. Lu and the author provided in \cite{LW19b} a realization of $\tUi$ via Nakajima-Keller-Scherotzke quiver varieties, generalizing F.~ Qin's approach for quantum groups \cite{Qin16}.

\item
\emph{$\imath$Categorification}

There has been a KLR type categorification of one family of modified $\imath$quantum groups of type AIII by Bao, Shan, Webster and the author \cite{BSWW}. The categorification of the split rank one $\imath$quantum group (see Example~\ref{ex:rank1}) will be a fundamental new step, allowing the $\imath$categorification to move forward. The $\imath$categorification shall have applications to modular representation theory. 

\item
\emph{$\imath$Hall algebra}

So far, the $\imath$Hall algebras can only realize the \emph{quasi-split} $\imath$quantum groups; see \cite{LW20a}. It is desirable to extend the $\imath$Hall algebras to a greater generality allowing Satake diagrams with black nodes, and also to understand categorically the embedding $\tUi \rightarrow \tU$ as well as the coideal structure ($\Delta:\tUi \rightarrow \tUi \otimes \tU$). 

\item
\emph{Representations of affine $\imath$quantum groups}

There has been numerous results in finite-dimensional representations of affine quantum groups and connections to other areas, by V.~Chari and many others. One hopes the Drinfeld type presentation of affine $\imath$quantum groups (see \cite{LW20b, Zh21}) can stimulate the development of their finite-dimensional representations. 

\item
\emph{$\imath$Quantum groups at roots of 1}

Building on Lusztig's constructions for quantum groups at roots of 1, Bao and Sale \cite{BS21} have taken a first step in formulating small quantum symmetric pairs. More can be expected in this direction in light of Lusztig's program. 
\end{enumerate}

It is hoped that $\imath$quantum groups may find more applications in math physics, geometric and modular representation theory, quantum topology, and algebraic combinatorics. 


Just as generalizing the study of compact or complex Lie groups to real Lie groups and symmetric spaces, we hope to have convinced the reader that it is a good thing to do to generalize various fundamental constructions from quantum groups to $\imath$quantum groups. 

It is time for the reader to come up with his own favorite item (+) in the list of highlights for quantum groups in the Introduction, and supply its missing $\imath$-generalization!

\vspace{2mm}

{\bf Acknowledgments.} 
This survey represents part of what the author has learned about $\imath$quantum groups from his friends, many collaborators and students, to whom he is very grateful;
 special thanks go to Huanchen Bao, Stefan Kolb, Ming Lu, and Weinan Zhang.
The dramatic development and rapid expansion in this research area in recent years are made possible largely due to their ideas, insights, generosities, and tireless efforts. We thank Huanchen Bao, Ming Lu, and Hideya Watanabe for careful reading of the paper and helpful comments. 
This work was partially supported by the NSF grant DMS-2001351.


\end{document}